\documentclass[11pt]{amsart}
\usepackage{amsmath}
\usepackage{amssymb}
\usepackage{amsfonts}
\usepackage{mathrsfs}
\usepackage{color}
\usepackage{graphics}
\usepackage[dvips]{graphicx}
\usepackage{wrapfig}
\usepackage{upgreek}
\usepackage[colorinlistoftodos, dvistyle, textsize=tiny]{todonotes} 
\usepackage[mathscr]{euscript}

\newtheorem{theorem}{Theorem}[section]

\newtheorem{definition}{Definition}[section]
\newtheorem{remark}[theorem]{Remark}

\newcommand{\p}{\partial}
\newcommand\R{\mathbb R}

\renewcommand\det{\text{det}}

\newcommand\ol{\overline}

\renewcommand\phi{\varphi}

\newcommand\cal{\mathcal}

\newcommand{\U}{\mathcal U}
\newcommand{\V}{\mathcal V}
\renewcommand{\S}{\mathbb{S}}

\renewcommand\phi\varphi 
\numberwithin{equation}{section}
\newcommand\M{\mathscr M}

\numberwithin{equation}{section}

\newcommand\qf{\text{\textnormal{II}}}
\newcommand\T{\mathbb T}
\newcommand\X{\mathscr X}

\DeclareMathAlphabet{\mathpzc}{OT1}{pzc}{m}{it}
\usepackage{sidecap}
\usepackage{subfig}

\title[Blaschke's rolling ball theorem and the Trudinger-Wang monotone bending]
{Blaschke's rolling ball theorem and the Trudinger-Wang monotone bending }
\author[Aram L. Karakhanyan ]{Aram L. Karakhanyan }
\address{Aram L. Karakhanyan\\
School of  Mathematics\\
The University of Edinburgh\\
Mayfield Road, EH9 3JZ\\
Edinburgh, UK}
\email{aram.karakhanyan@ed.ac.uk}

\thanks{ 
\\
 }

\begin{document}

\maketitle

\baselineskip=15pt

\begin{abstract}
We revisit the classical rolling ball theorem of  Blaschke for 
convex surfaces with positive curvature  and show that  
it is linked to  another inclusion principle in the 
optimal mass transportation theory due to  
Trudinger and Wang. We also discuss an application to reflector antennae design problem.
\end{abstract}

\section{Introduction}

In this note we give two applications of an inclusion principle known as    
the rolling ball Theorem of Blaschke. 
Let $\M$ and $\M'$ be two hypersurfaces in $\R^d$. We say that 
$\M$ and $\M'$ are {\it internally tangent at} $x\in \M$ if 
they are tangent at $x$ and have the same outward normal. 
Denote by $\qf_x\M$ the second fundamental form of $\M$ at $x$ and 
let $n(x)$ be the outward unit normal at $x$. Then we have 

\begin{theorem}\label{th-1}
Suppose $\M$ and $\M'$ are smooth convex surfaces with strictly 
positive scalar curvature  such that $\qf_x\M\ge \qf_{x'}\M'$ for all 
$x\in\M, x'\in \M'$ such that $n(x)=n'(x')$. If $\M$ and $\M'$ are internally tangent at one point 
then $\M$ is contained in the convex region bounded by $\M'$.
\end{theorem}
\medskip 
W. Blaschke \cite{Blaschke}, pp. 114-117 proved Theorem \ref{th-1}  for closed curves in $\R^2$. 
D. Koutroufiotis \cite{Kou} generalized Blaschke's theorem  for complete curves in $\R^2$ 
and complete surfaces in $\R^3$. Later J. Rauch \cite{Rauch}, by using Blaschke's techniques, proved 
this result  for compact surfaces in $\R^d$ and J.A. Delgado \cite{Del} for  complete surfaces. Finally J. N. Brooks and  J. B. Strantzen generalized  Blaschke's theorem 
 for  non-smooth convex sets showing that the local inclusion implies global inclusion \cite{Brooks}.

\smallskip

Observe that if  $\M$ and $\M'$ are  internally tangent at $x$, then a necessary condition for $\M$ to be inside 
 $\M'$  near $x$ is  
\begin{equation}
\qf_x(v)\ge \qf'_{x}(v)\quad \text{for all} \ v \in \T_x\M\cong\T_{x}\M'.
\end{equation}
The tangent planes are parallel because $\M$ and $\M'$ are internally tangent at $x$.
Therefore Theorem \ref{th-1} says that if  for all $x\in \M, x'\in\M', x\not=x'$ with coinciding 
normals  $n'(x')=n(x)$
such that  after translating $\M$ by $x-x'$ 
we have that the translated surface $\widetilde \M$ is locally inside $\M'$ then
$\M$ is globally inside  $\M'$.  In other words, the local inclusion implies global inclusion
or $M$ rolls freely inside $\M'$.

\smallskip 

Our aim is to apply Theorem \ref{th-1} to  optimal transportation theory and 
reflector antennae design problems. More specifically, for a smooth cost function 
$c:\R^d\times \R^d\to \R$ (subject to some standard conditions including the weak $\bf A3$ condition) 
and a pair of bounded smooth convex domains $\U, \V\subset \R^d$ such that 
$\U$ is $c-$convex with respect to $\V$ (see Definition \ref{def-long} below), we would like to take 
$\M=\U$ to be the reference domain and 
$$\M'=\mathcal N:=\{x\in \R^d \ s.t. \ c(x, y_0)=c(x, y_1)+a\}, \quad y_1, y_2\in \V$$ 
for some  constant $a$. 
Then $\M'$ is the boundary of sub-level set of the cost function $c$.
We prove that if $\p\U$ is locally inside $\mathcal N$ in above sense then 
$\p\U$ is globally inside $\mathcal N$ provided that the  sets $\mathcal N$ are convex for all $y_1, y_2, a$. 
The precise result is formulated in Theorem \ref{th-2} below and applications in 
Section \ref{sec-app}.

A local inclusion principle for $\p\U$ and $\mathcal N$ is proved by Neil Trudinger and Xu-Jia Wang in \cite{TW-2}, see 
the inequality (2.23) there. It is then used to show that 
under the $\bf A3$ condition a local support function is also global, see \cite{TW-2} page 411. 
The proof is based on a monotone bending argument that gives yet another geometric interpretation of the
$\bf A3$ condition.


\section{Preliminaries}

\subsection{Optimal transportation} In order to formulate the result in the context of  optimal transportation theory we 
need some standard definitions.
Let $c:\R^d\times\R^d\to \R$ be a cost function such that $c\in C^4(\R^d\times\R^d)$ and 
$\U, \V\subset \R^d$.
\begin{definition}\label{def-long}
\ 

\begin{itemize}
\item Let $u:\U\to \R$ be a continuous function.
A $c-$support function of $u$ at $x_0\in \U$ is $\phi_{x_0}=c(x, y_0)+a_0, y_0\in \R^d$ such that the following two conditions hold
\[\begin{array}{lll}
u(x_0)&=&\phi_{x_0}(x_0),\\
u(x)&\ge& \phi_{x_0}(x), x\in \U.
\end{array}\]
\item If $u$ has $c-$support at every $x_0\in \U$ then we say that 
$u$ is $c-$convex in $\U.$
\item $c-$segment with respect to a point $y_0\in \R^d$ is the set 
$$\{x\in \R^d\ s.t. \ c_y(x, y_0)=\text{line segment}\}.$$
One may take in the above definition $\{x\in \R^d\ s.t. c_y(x, y_0)=tp_1+(1-t)p_0\}$
with $t\in[0, 1]$ and $p_0, p_1$ being two points in $\R^d$.
\item 
We say that $\U$ is $c-$convex with respect to $\V\subset \R^d$ if the image of the set 
$\U$ under the mapping $c_y(\cdot, y)$ denoted by 
$c_y(\U, y)$ is convex set for all $y\in \V.$
Equivalently, $\U$ is $c-$convex with respect to $\V$ if for any pair of points $x_1, x_2\in \U$
there is $y_0\in \V$ such that there is a  $c-$segment with respect to $y_0$ joining $x_1$ with $x_2$
and lying in $\U$.
\end{itemize}
\end{definition}

\begin{definition} 
Let $u$ be a $c-$convex function then the sub-level set of $u$ at $x_0\in \U$ is 
\begin{equation}
S_{h, u}(x_0)=\{x\in \R^d\ s.t. u(x)<c(x, y_0)+[u(x_0)-c(x_0, y_0)]+h\}
\end{equation}
for some constant $h$. 
\end{definition}

Equivalently, $S_{h, u}(x_0)=\{x\in \U\ s.t.u(x)< \phi_{x_0}(x)+h\}$ where 
$\phi_{x_0}$ is the $c-$support function of $u$ at $x_0\in \U$, see Definition \ref{def-long}.

Observe that in the previous definition on may take $u(x)=c(x, y_1)$ for some 
fixed $y_1\not=y_0$. 

Next we recall Kantorovich's formulation of optimal transport problem, see \cite{Urbas, Vil}:
Let $f:\U\to\R, g:\V\to \R$ be two nonnegative integrable functions satisfying the mass balance condition 
\begin{equation*}
\int_\U f(x)dx=\int_\V g(y)dy.
\end{equation*} 
Then one wishes to minimize 
\begin{equation}
\int_\U u(x)f(x)+\int_\V v(y)g(y)dy\to \min
\end{equation}
among all pairs of functions $u:\U\to \R, v:\V\to \R$  such that 
$u(x)+v(y)\ge c(x,y).$

It is well-known that a minimizing pair $(u, v)$ exists \cite{Urbas, Vil} and  formally the potential $u$ solves the equation 
\begin{equation}
\det(u_{ij}-A_{ij}(x, Du))=|\det c_{x_i, y_j}|\frac{f(x)}{(g\circ y)(x)}. 
\end{equation}
Here $A_{ij}(x, p)=c_{x_ix_j}(x, y(x,p))$ where $y(x, p)$ is determined from 
$D_x(c(x, y(x, p)))=p$.

Assume that $c$ satisfies the following conditions:
\begin{itemize}
\item[{$\bf A1$}] For all $x, p\in \R^d$ there is unique $y=y(x, p)\in \R^d$ such that $\partial_x c(x, y)=p$
and for any $y, q\in \R^d$ there is unique $x=x(y, q)$ such that $\partial_yc(x, y)=q$.
\item[{$\bf A2$}] For all $x, y\in \R^d$ $\det c_{x_i, y_j}(x, y)\not=0$.
\item[{$\bf A3$}] For $x, p\in \R^d$ there is a positive constant $c_0>0$ such that 
\begin{equation}\label{MTW}
A_{ij,kl}(x, p)\xi_i \xi_j\eta_k\eta_l\ge c_0|\xi|^2|\eta|^2\quad \forall \xi, \eta\in \R^d, \xi\perp\eta.
\end{equation}
\end{itemize}
{\bf A3}  is the Ma-Trudinger-Wang condition \cite{MTW}.

J.Liu proved that if {\bf A1-A3} hold then $S_{h, u}(x_0)$ is $c-$convex
with respect to $y_0$  \cite{Liu}. There are cost functions satisfying the weak {\bf A3}
\begin{equation}\label{MTW-0}
A_{ij,kl}(x, p)\xi_i \xi_j\eta_k\eta_l\ge 0\quad \forall \xi, \eta\in \R^d, \xi\perp\eta.
\end{equation}
i.e. when $c_0=0$ in \eqref{MTW}, such that 
the corresponding sub-level sets are convex in classical sense, see Section \ref{sec-app}.

We also remark that the condition 
{\bf A3} is equivalent to 
\begin{equation}\label{monot}
\frac{d^2}{dt^2}c_{ij}(x, y(x, p_t))\xi_i\xi_j\ge c_0|p_1-p_0|^2
\end{equation}
where $x$ is fixed, $c_x(x, y(x, p_t))=tp_1+(1-t)p_0, t\in[0, 1]$ 
$c_x(x, y)=p_1, c_x(x, y_0)=p_0$ (this determines the so-called $c^*-$segment
with respect to fixed $x$), see \cite{TW-2}.

\subsection{Shape operator}
If $\M$ is a surface with positive sectional curvature then by Sacksteder's theorem \cite{Sack}
$\M$ is convex.
For  $x\in \M$, let $n(x)$ be the unit outward normal 
at $x$ ($n(x)$ points outside of the convex body bounded by $\M$). 
The Gauss map $x\to n(x)$ is a diffeomorphism of $\M$ onto $\S^{d}$ \cite{Wu}, 
where $\S^d$ is the unit sphere in $\R^d$. The inverse map  $n^{-1}$ gives a parametrization of $\M$ by $\S^d$. 
 If $\M'$ is another smooth convex  surface, and $w\in  \S^d$, then $n^{-1}(w)$ and $(n')^{-1}(w)$ are the points
 on $\M$ and $\M'$ with equal outward normals.

 Let $F:\Omega\to \R^m$ be a smooth map on a set $\Omega\subset \R^d$
 and $v=(v_1, \dots, v_d)\in\R^d$ then 
 $$\partial_v F(x)=\sum_{i=1}^{d}v_i\frac{\partial F(x)}{\partial x_i}$$
 is the directional derivative operator.

We view the tangent space as a linear subspace of $\R^d$ consisting of tangential directions. Then 
the tangent space $\T_x\M$ is the set of vectors perpendicular to $n(x)$.

\smallskip 

Next we introduce the Weingarten map in order to define the second fundamental form.
The Weingarten map $W_x:\T_x\to \T_x$ is defined by $W_x(v)=\p_v n(x)$.
$\T_x$ is an inner product space (induced by the inner product in $\R^d$).
Then $W_x$ is self-adjoint operator on $\T_x$ and the eigenvalues of
$W_x$ are the principal curvatures at $x$.

\begin{remark}\label{rem-on}
Observe that since $W_x$ is self-adjoint 
and $\T_x$ is finite dimensional then there exists an orthonormal basis of
$\T_x$ consisting of eigenvectors of $W_x$.
\end{remark}

\begin{definition}
The second fundamental form is defined
as $\qf_x(v, w)=W_x(v)\cdot w$. When 
$v=w$ we denote $\qf_x(v)$.
\end{definition}

From definition it follows that if $\M$ is parametrized by $r=r(u)$ and $x=r(u_0)$ then  
\begin{equation}\label{qf}
\qf_{x}(v)=-\p_v^2 r \cdot n(x), \quad v\in \T_x
\end{equation}
which readily follows from the differentiation of $n\cdot \p_v r=0$.
\section{Main result}
\begin{theorem}\label{th-2}
Let $y_1, y_2\in \ol{\V}$ and $\cal N(y_1, y_2, a)=\{x\in \R^d : c(x, y_0)=c(x, y_1)+a\}$ for some $a\in \R$
where $c$ satisfies $\bf A1, A2$ and weak $\bf A3$, see \eqref{MTW-0}.
Assume that $\mathcal N$ is convex for all $y_1, y_2, a$ and $\U$ is convex domain with 
smooth boundary such that $\U$ is $c-$convex with respect to $\V$, see Definition \ref{def-long}.
If $\cal N$ and $\p\U$ are internally tangent at some point $z_0$ then $\U$ is inside $\cal N$.
\end{theorem}
Using the terminology of Blaschke's theorem it follows that under the conditions of Theorem \ref{th-2} $\U$ rolls freely inside $\mathcal N$.
Observe that the $c-$convexity of sub-level sets is known under stronger condition $\bf A3$ \cite{Liu}. In the next 
section we give an example of cost function $c$ satisfying weaker form of $\bf A3$ \eqref{MTW-0} but such that  $\mathcal N$ is convex for all $y_1, y_2, a$.
Proof to follow is inspired in \cite{TW-2}.

\begin{proof}

{\bf Step 1:} ({\it Parametrizations})

\smallskip 
To apply Theorem \ref{th-1} we take $\M=\U$ and $\M'=\cal N$ and assume that $\U$ and 
$\mathcal N$ are internally tangent at $z_0$. Assume that at $x_0'\in \cal N$ and $x_0\in \p\U$ 
$\cal N$ and $\p\U$
have the same outward normal, see Figure 1. 

In what follows we use the following radial parametrizations:
\[
\begin{array}{lll}
\p\U& \quad R(\zeta), \quad \zeta\in D_\U,\\ 
\cal N& \quad \X(\omega), \quad \omega\in  D_{\cal N},\\
\p(c_y(\p\U, y_0))& \quad \rho(\zeta)=c_y(R(\zeta), y_0).
\end{array}
\]
Here $D_\U$ and $D_{\mathcal N}$ are the domains of corresponding parameters.
Moreover, there are $\bar\omega\in D_{\mathcal N}$ and $\bar \zeta\in D_{\U}$ such that
\begin{equation}
 x_0':=\X(\bar \omega)\in \cal N \quad \text{and}\quad  x_0:=R(\bar \zeta)\in \p\U.
\end{equation}
\smallskip

From now on $\bar \zeta$ and $\bar\omega$ are fixed. 
Let $\bar n(\bar \zeta)$ denote the outward normal of the image $c_y(\U, y_0)$ 
at the point $\rho(\bar \zeta)$. We have 
\begin{equation}\label{normal}
\bar n^m(\bar\zeta)=c_{y_m, x_i}(R(\bar \zeta), y_0)n^i(\bar \zeta).
\end{equation}
Observe that by assumption the constant matrix $\mu=[c_{y_m, x_i}(R(\bar \zeta), y_0)]^{-1}$ 
has non-trivial determinant,  see {\bf A2}. Furthermore, the set 
$\mu c_y(\U, y_0)=\{\mu x\ s.t. \ x\in c_y(\U, y_0)\}$ is again convex because 
for any two points $q_1=\mu z_1, q_2=\mu z_2$ such that 
$q_1, q_2\in \mu c_y(\U, y_0)$ and $z_1, z_2\in c_y(\U, y_0)$ we have 
$$\mu c_y(\U, y_0)\ni \mu(\theta z_1+(1-\theta)z_2)=\theta\mu z_1+(1-\theta)\mu z_2= \theta q_1+(1-\theta)q_2 $$
for all $\theta\in[0, 1]$.

\begin{figure}\label{fig-1}
\includegraphics[scale=0.40]{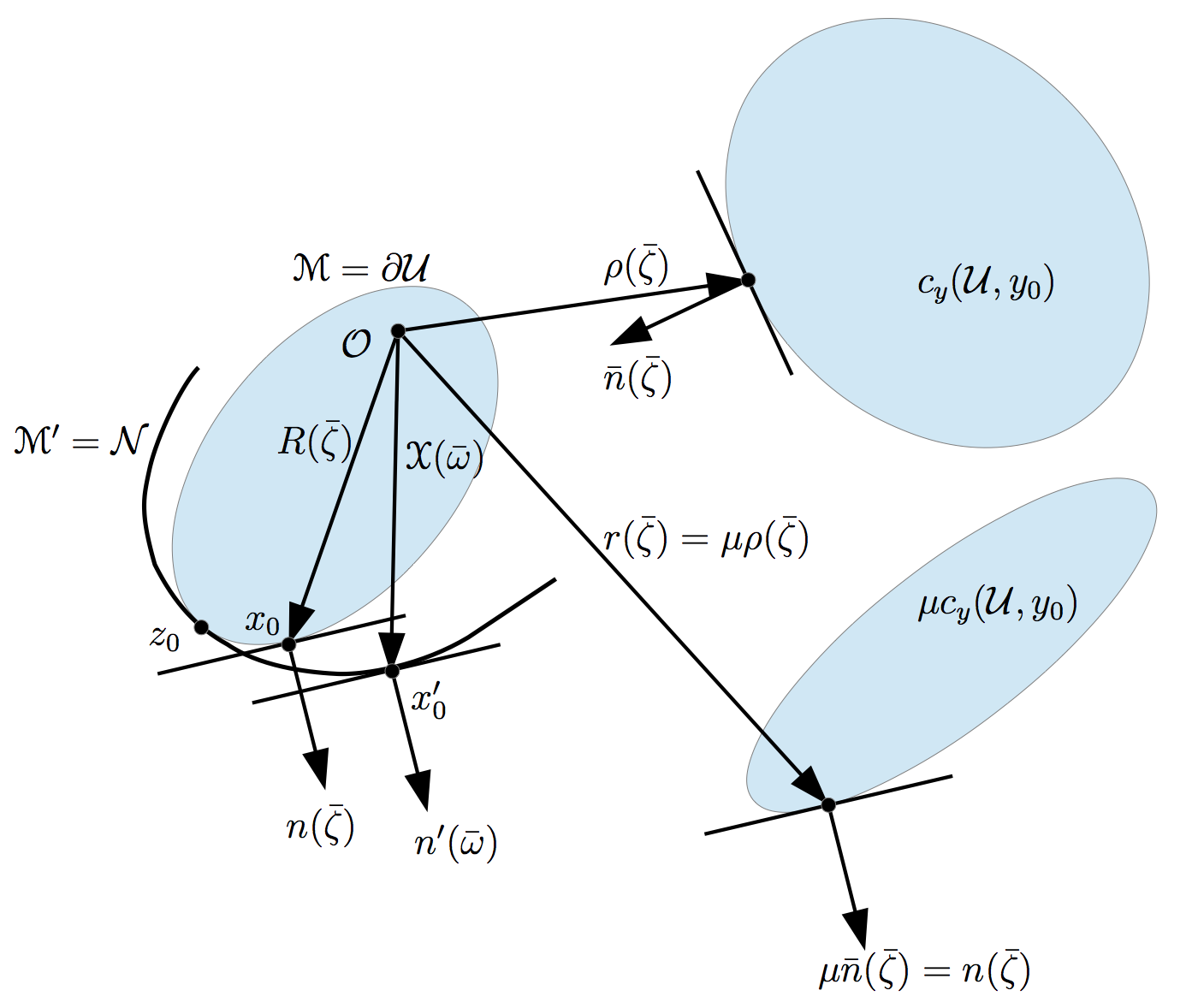}
\caption{Schematic view to parametrizations of $\p\U, \mathcal N, \p (c_y(\U, y_0))$ and $ \mu \p( c_y(\U, y_0))$.}
\end{figure}

\medskip 
{\bf Step 2:} ({\it Computing the second fundamental form of $\X$})

Next, we introduce the vectorfield $r=r(\zeta), \zeta\in D_\U$  such that 
\begin{equation}
r(\zeta)=\mu\rho(\zeta)=\mu c_y(R(\zeta), y_0).
\end{equation} 
We compute the first and second derivatives 
\begin{eqnarray}\label{norm-1}
r_{\zeta_s}^m&:=&r_s^m=\mu_{\alpha\beta} c_{y_\beta, x_i} R^i_s,\\
r_{st}^m&=& \mu_{\alpha\beta}\left[c_{y_\beta, x_ix_j}R^i_sR^j_t+c_{y_\beta, x_i}R^i_{st}\right].
\end{eqnarray}

From \eqref{norm-1} and \eqref{normal} we see that 
at $r(\bar\zeta)$ the normal is 
\begin{equation}\label{norm-3}
n(\bar \zeta)=\mu\bar n(\bar\zeta).
\end{equation}

Take $p_t=(1-t)p_0+tp_1, t\in[0,1]$ and 
\begin{equation}\label{p-t}
p_t=c_x(x_0, y(x_0', p_t)), 
\end{equation}
then $y_t:= y(x_0', p_t)$ defines the $c-$segment joining $y_0$ and $y_1$, see {\bf A2}.
In particular, one has 
\begin{eqnarray}\label{y}
p_1^i-p_0^i&=&c_{x_i, y_m}(x_0, y(x_0', p_t))\frac{d}{dt}y^m(x_0', p_t)\\\nonumber
&=&\frac{d}{dt}y^m(x_0', p_t)c_{y_m, x_i}(x_0, y(x_0', p_t)).
\end{eqnarray}

Let $\X^t(\omega)$ be the parametrization of $\cal N(t)=\{x\in \U : c(x, y_0)=c(x, y_t)+a\}$ (recall that 
$\mathcal N(t)$ is convex as the boundary of sub-level set).
We can choose $a=a(t)$ so that 
all $\cal N(t)$ pass through the point $x_0'$, in other words there is $\bar\omega^t$ such that 
$\X^t(\bar\omega^t)=x_0'$.
Moreover, by \eqref{p-t} it follows that  
\begin{eqnarray}
c_{x_i}(\X^t(\bar\omega^t), y_0)-c_{x_i}(\X^t(\bar\omega^t), y_t)&=&c_{x_i}(x_0', y_0)-c_{x_i}(x_0', y_t)\\\nonumber
&=&p_0^i-p_t^i\\\nonumber
&=&t(p_0^i-p^i_1).
\end{eqnarray}

After fixing $t$ and differentiating the identity $c(\X^t(\omega), y_0)=c(\X^t(\omega), y_t)+a(t)$ in $\omega$ we get
\begin{eqnarray}\nonumber
\left[c_{x_i}(\X^t, y_0)-c_{x_i}(\X^t, y_t)\right]\X_{\omega_k}^{i, t}=0,\\\label{long}
\left[c_{x_ix_j}(\X^t, y_0)-c_{x_ix_j}(\X^t, y_t)\right]\X^{j, t}_{\omega_l}\X_{\omega_k}^{i, t}+
\left[c_{x_i}(\X^t, y_0)-c_{x_i}(\X^t, y_t)\right]\X_{\omega_k\omega_l}^{i, t}=0.
\end{eqnarray}
Thus the normals of $\cal N(t)$ at $x_0'$ are collinear to $p_1-p_0$ for all $t\in[0,1]$, that is  

\begin{equation}\label{nnn}
n(x_0)=n'(x_0')=\frac{p_1-p_0}{|p_1-p_0|}, \quad \mu \bar n=n \ \ (\text{recall}\ \ \eqref{norm-3}).
\end{equation}

Hence we can rewrite \eqref{long} as follows
\begin{eqnarray}
\left[(c_{x_ix_j}(\X^t, y_0)-c_{x_ix_j}(\X^t, y_t)\right]\X^{j, t}_{\omega_l}\X_{\omega_k}^{i, t}=
-t(p_0^i-p^i_1)\X_{\omega_k\omega_l}^{i, t}.
\end{eqnarray}
Keeping $\X^t(\bar \omega^t)=x_0'$ fixed for all $t\in[0,1]$, dividing both sides of the last identity 
by $t$ and then sending $t\to 0$ we obtain 
\begin{eqnarray}
-\left[y'(x_0', p_0)c_{y, x_ix_j}(x_0', y_0)\right]\X^{j, t=0}_{\omega_l}\X_{\omega_k}^{i, t=0}=
-(p_0^i-p^i_1)\X_{\omega_k\omega_l}^{i, t=0}.
\end{eqnarray}
On the other hand from \eqref{y} we see that $\frac d{dt}y(x_0', p_t)\big|_{t=0}=(p_1-p_0)\mu$. Thus substituting 
this into the last equality we obtain 
\begin{eqnarray}\label{blya}
\left[(p_1-p_0)\mu c_{y, x_ix_j}(x_0', y_0)\right]\X^{j, t=0}_{\omega_l}\X_{\omega_k}^{i, t=0}&=&
(p_0^i-p^i_1)\X_{\omega_k\omega_l}^{i, t=0}\\\nonumber
&=&-(p_1^i-p^i_0)\X_{\omega_k\omega_l}^{i, t=0}
\end{eqnarray}
or equivalently 
\begin{eqnarray}\label{blya-1}
\left[ n^\alpha \mu_{\alpha\beta} c_{y_{\beta}, x_ix_j}(x_0', y_0)\right]\X^{j, t=0}_{\omega_l}\X_{\omega_k}^{i, t=0}=
-n^i\X_{\omega_k\omega_l}^{i, t=0}\\\nonumber
\end{eqnarray}
if we utilize \eqref{nnn}.

\medskip 
{\bf Step 3:} ({\it Monotone bending})

By Remark \ref{rem-on} we assume that 
$\T_{x_0}\partial \U$ and $\T_{x'_0}\mathcal N(t=0)$ have the same local coordinate system (by reparametrizing $\mathcal N(t=0)$ if necessary). 
From convexity of $\mu c_y(\U, y_0)$ boundary of which is parametrized by 
$r$ we have

\begin{eqnarray}\label{blya-3}
0&\ge& r_{st}^\alpha n^\alpha=\mu \rho_{st}  n =\mu_{\alpha\beta}(c_{y_\beta, x_ix_j}R^i_sR^j_t+c_{y_\beta, x_i}R^i_{st})  n^\alpha\\\nonumber
&=&\mu_{\alpha\beta} c_{y_\beta, x_ix_j}R^i_sR^j_t  n^\alpha+R^i_{st}  n^i\\\nonumber
&\overset{\eqref{blya-1}}{=}&-n^i \X_{\omega_k\omega_l}^{i, t=0}+ R^i_{st} n^i. 
\end{eqnarray}
Now \eqref{monot} yields that at $x_0'$ 
\begin{eqnarray}
n^i \X_{\omega_k\omega_l}^{i, t}\ge n^i \X_{\omega_k\omega_l}^{i, t=0} \overset{\eqref{blya-3}}{\ge}  R^i_{st}  n^i. 
\end{eqnarray}
Recalling \eqref{qf} we finally obtain the required inequality
$$\qf_{x_0'} \mathcal N\leq \qf_{x_0}\p\U.$$
The proof is now complete.
\end{proof}
Note that  weak {\bf A3} (i.e. when $c_0=0$ in \eqref{monot}) is enough for the monotonicity to conclude the inequality
$n^i \X_{\omega_k\omega_l}^{i, t}\ge n^i \X_{\omega_k\omega_l}^{i, t=0}.$


\section{Applications}\label{sec-app}
\subsection{Convex sub-level sets} 
There is a wide class of cost functions for which the set 
$\cal N$ is convex. Observe that  $c(x, y)=\frac1p{|x-y|^p}$
satisfies {\bf A3} for $-2<p<1$ and weak {\bf A3} if $p=\pm2$ \cite{MTW}.

It is useful to note that if  $\Omega_\psi=\{x\in \R^d\ s.t.\ \psi(x)<0\}$ for some 
smooth function $\psi:\R^d\to \R$ such that $\Omega_\psi\not=\emptyset$
then 
\begin{equation}\label{hessian}
\partial^2 \psi(x)\tau(x)\cdot \tau(x)\ge 0, \quad \forall \tau(x)\in \T_x
\end{equation}
is a necessary and sufficient condition for $\Omega_\psi$
to be convex provided that $\frac{\nabla \psi}{|\nabla \psi|}$ is directed towards positive $\psi$.

Using this, one can check that the boundaries of sub-level sets (e.g. $p=-2$)
$$\frac1{|x-y_2|^2}-\frac1{|x-y_1|^2}-a=0$$
are convex. Here $\psi(x):=|x-y_2|^{-2}-|x-y_1|^{-2}-a<0$ defines the sub-level set.
Indeed, using \eqref{hessian} we see that 
for any unit vector $\tau$ perpendicular to the vector 
$$\nabla\psi(x)=-\frac{2(x-y_2)}{|x-y_2|^4}+\frac{2(x-y_1)}{|x-y_1|^4}$$
we get 
\begin{eqnarray*}
\partial^2_{\tau\tau}\psi&=&-\frac2{|x-y_2|^4}\left[1-\frac{4((x-y_2)\cdot\tau)^2}{|x-y_2|^2}\right]
+ \frac2{|x-y_1|^4}\left[1-\frac{4((x-y_1)\cdot\tau)^2}{|x-y_1|^2}\right]\\\nonumber
&\overset{ \nabla\psi\cdot\tau=0 }{=}&\frac2{|x-y_1|^4}-\frac2{|x-y_2|^4}-\frac{8((x-y_1)\cdot\tau)^2}{|x-y_1|^6}+\frac{8((x-y_1)\cdot\tau)^2}{|x-y_2|^6}\frac{|x-y_2|^8}{|x-y_1|^8}\\\nonumber
&=&\frac2{|x-y_1|^4}-\frac2{|x-y_2|^4}+ \frac{8((x-y_1)\cdot\tau)^2}{|x-y_1|^6}\left[\frac{|x-y_2|^2}{|x-y_1|^2}-1\right]
\\\nonumber
&=&\frac2{|x-y_1|^4}-2\left(\frac1{|x-y_1|^2}+a\right)^2+  \frac{8((x-y_1)\cdot\tau)^2}{|x-y_1|^6}
\left[ 
\frac{\frac1{|x-y_1|^2}}{\frac1{|x-y_1|^2}+a}-1 
\right]\\\nonumber
&=&\frac2{|x-y_1|^4}-2\left(\frac1{|x-y_1|^2}+a\right)^2-  \frac{8((x-y_1)\cdot\tau)^2}{|x-y_1|^6} 
\frac{a}{\frac1{|x-y_1|^2}+a}.
\end{eqnarray*}
Altogether, we infer that 
\[
\partial^2_{\tau\tau}\psi
\begin{array}{lll}
\le 0\quad  \text{if}\ a\ge0.
\end{array}
\]
If $a<0$ then we can swap $y_1$ and $y_2$ to conclude that $\partial^2_{\tau\tau}\psi>0$.
Some examples of sub-level sets of inverse quadratic cost function are illustrated in Figure 2.

\begin{figure}
\includegraphics[scale=0.4]{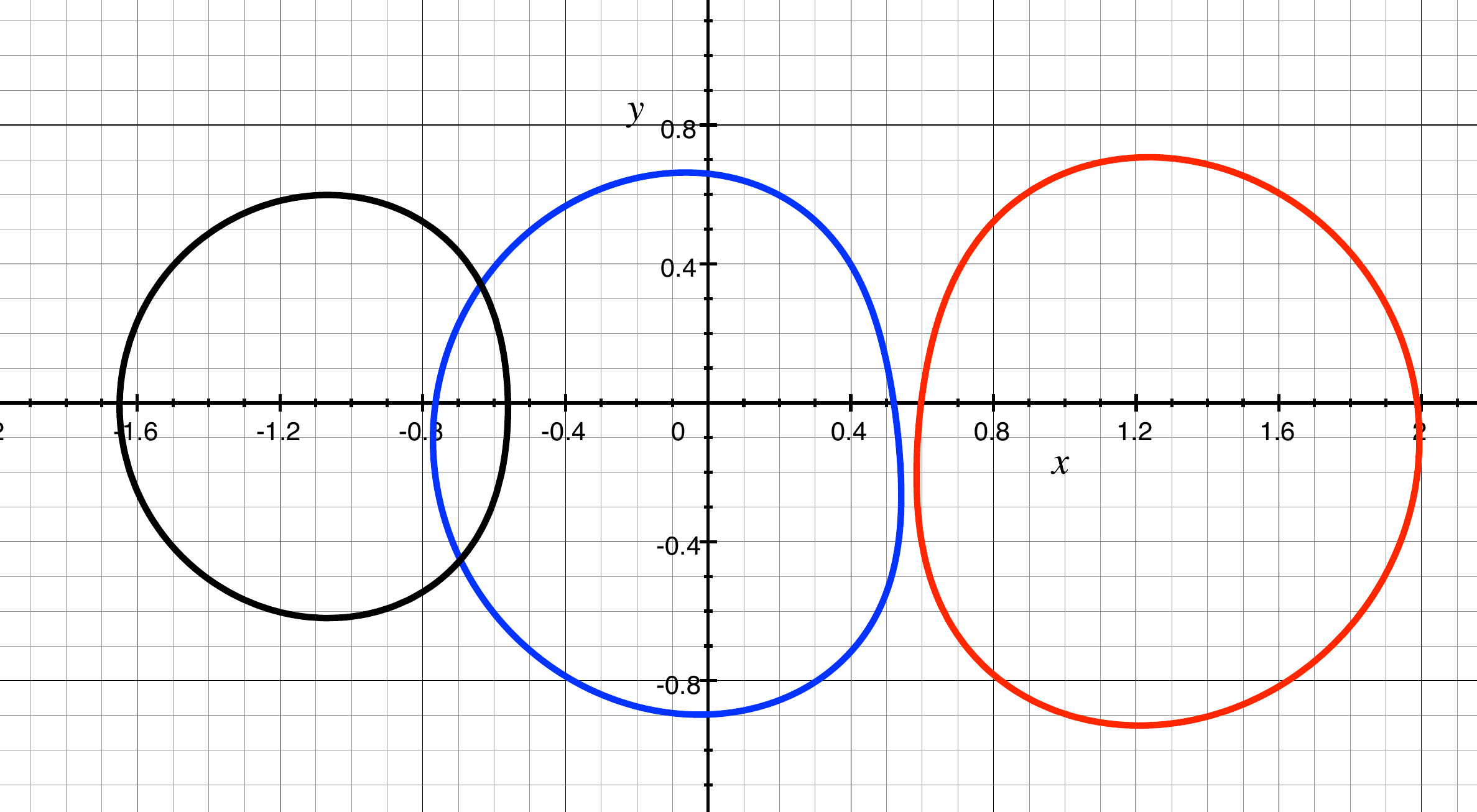}
\caption{From left to right: $a=-2, y_1=(-10^{-3}, 0), y_2=(-1, -10^{-2})$; 
$a=1, y_1=(-10^{-1}, -10^{-1}), y_2=(1, 10^{-2})$; $a=-1, y_1=(-10^{-4}, 0), y_2=(1.1, -10^{-1}).$}
\end{figure}
\subsection{Antenna design problems}
In parallel reflector problem 
\cite{K} one deals with the paraboloids of revolution 
\begin{equation}\label{parab-gen-f}
 P(x,\sigma, Z)=\frac\sigma2+Z^{n+1}-\frac1{2\sigma}|x-z|^2
\end{equation}
which play the role of support functions. 
Here the point $Z=(z, Z^{n+1})\in \R^{n+1}$ is the focus of 
the paraboloid such that $\psi(z, Z^{n+1})=0$ for 
some smooth function $\psi$ satisfying some structural conditions and 
$\sigma$ is a constant. If $P_1$ is internally 
tangent to $P_2$ at $z_0$ and $\qf_{z_0} P_1\ge \qf_{z_0} P_2$ then 
$P_1$ is inside  $P_2$, see Lemma 8.1 \cite{K}.  
This again follows from Blaschke's theorem. Indeed, we have that at the 
points $x$ and $x'$ corresponding to coinciding outward normals 
$$\qf_x P_1=\frac1{\sqrt{1+|DP_1(x)|^2} } \frac{1}{\sigma_1}\delta_{ij}$$ 
and $$\qf_{x'} P_2=\frac1{\sqrt{1+|DP_2(x')|^2}} \frac{1}{\sigma_2}\delta_{ij}.$$
Furthermore $DP_1(x)=DP_2(x')$ and hence 
\begin{equation}\label{koren}
\sqrt{1+|DP_1(x)|^2}=\sqrt{1+|DP_2(x')|^2}.
\end{equation} 
From $\qf_{z_0} P_1\ge \qf_{z_0} P_2$ we infer that 
\begin{equation}\label{sigma}
\frac1{\sigma_1}\ge \frac1{\sigma_2}.
\end{equation}
Consequently 
\eqref{sigma} and \eqref{koren} imply that 
$$\qf_{x} P_1\ge \qf_{x'} P_2.$$

\subsection{Another inclusion principle}

There are various inclusion principles in geometry,  we want to mention the following elementary one due to 
Nitsche \cite{Nit}: Each continuous closed curve of length $L$ in Euclidean 3-space is contained in a 
closed ball of radius $R<L/4.$  Equality holds only for a "needle", i.e., a segment of length
$L/2$ gone through twice,  in opposite directions. Later J. Spruck generalized this result 
for compact Riemannian manifold $\M$ of dimension $n\ge 3$ as follows: if the  sectional curvatures 
$K(\sigma)\ge1/c^2$ for all tangent plane sections $\sigma$ then $\M$ is contained in a ball of radius $R<\frac12\pi c$, and this bound is best possible. 
We remark here that there is a smooth surface   $S\subset \R^3$  such that the 
mean curvature $H\geq 1$ and the Gauss curvature $K \geq 1$
then the unit ball cannot be fit inside $S$, see \cite{Spruck}. Notice that $K$ is an intrinsic quantity and 
$H\ge 1$ implies that  $K\ge 1$.

\medskip


\end{document}